\documentclass[conference]{IEEEtran}

\newcommand{\pathtofigs}{./}

\IEEEoverridecommandlockouts
\usepackage{bm}
\usepackage{cite}
\usepackage{amsmath,amssymb,amsfonts,mathtools}
\usepackage[noend]{algorithmic}
\usepackage{graphicx}
\usepackage{textcomp}
\usepackage[usenames,dvipsnames,svgnames,table]{xcolor}
\usepackage{algorithm}
\usepackage{units}
\usepackage{url}
\usepackage{romannum}
\def\BibTeX{{\rm B\kern-.05em{\sc i\kern-.025em b}\kern-.08em
    T\kern-.1667em\lower.7ex\hbox{E}\kern-.125emX}}

\usepackage{tikz}
\usetikzlibrary{decorations.pathmorphing,
  patterns,shapes,arrows,positioning, backgrounds, intersections}
\tikzstyle{Nodes}=[circle, draw=black, fill=orange!20, line width=1.5pt, minimum size=15pt]
\tikzstyle{arrow} = [line width=1pt,->,>=stealth]
\tikzstyle{axis} = [line width=1pt,->,>=stealth]

\RequirePackage{pgfplots}
\usepgfplotslibrary{fillbetween}

\tikzstyle{block}=[
       draw=white,
       thick,
       minimum width=1cm,
       minimum height=1cm,
       inner sep=4pt,
       text width = 7em, 
       text centered,
       fill = white,
]

\usepackage{hyperref}
\hypersetup{
  unicode=false,          
  pdftoolbar=true,        
  pdfmenubar=true,        
  pdffitwindow=false,     
  pdfstartview={FitH},    
  pdftitle={notes},    
  pdfauthor={},                
  pdfsubject={gradgen},                   
  pdfcreator={LaTeX},                            
  pdfproducer={},      
  pdfkeywords={automatic differentiation, gradgen},      
  pdfnewwindow=true,      
  colorlinks=true,        
  linkcolor=blue!90,    
  citecolor=red,         
  filecolor=magenta,      
  urlcolor=NavyBlue           
}

\newcommand{\gradgen}{\texttt{gradgen}}
\newcommand{\casadi}{\texttt{CasADi}}
\newcommand{\R}{{\rm I\!R}}

\newcommand{\N}{{\rm I\!N}}
\newcommand{\tr}{{\intercal}}

\DeclareMathOperator*{\minimise}{\mathbf{Minimise}}
\DeclareMathOperator*{\subjto}{\mathbf{subject\ to:}}

\definecolor{myred}{rgb}{0.8,0.0,0.0}
\definecolor{mygreen}{rgb}{0.0,0.6,0.0}
\definecolor{myblue}{rgb}{0.0,0.0,0.8}

\usepackage{listings}
\lstdefinestyle{myRust}{
    language=c, 
    upquote=true,
    frame=single,    
    basicstyle=\small\ttfamily,
    backgroundcolor=\color{yellow!20},
    keywordstyle=[1]\color{NavyBlue}\bfseries,
    keywordstyle=[2]\color{purple},
    keywordstyle=[3]\color{orange},
    keywordstyle=[4]\color{Green},
    identifierstyle=,
    commentstyle=\usefont{T1}{pcr}{m}{sl}\color{blue}\small,
    stringstyle=\color{purple},
    showstringspaces=false,
    tabsize=5,              
    morekeywords={pub,fn,mut,let},
    morekeywords=[2]{f64,i32},
    morekeywords=[3]{vec},
    morekeywords=[4]{NX, NU, NPRED},
    otherkeywords={\&,!},
    morecomment=[l][\color{blue}]{...},
    numbers=none,
    firstnumber=1,
    numberstyle=\tiny\color{blue},
    stepnumber=1,
    xleftmargin=10pt,
    xrightmargin=10pt
}

\DeclareMathOperator{\child}{\mathbf{ch}}			
\DeclareMathOperator{\anc}{\mathbf{anc}}				
\DeclareMathOperator{\nodes}{\mathbf{nodes}}			

\newtheorem{theorem}{Theorem}
\newtheorem{proposition}[theorem]{Proposition}

\begin{document}

\title{Parallelisable computation of the gradient in nonlinear stochastic optimal control problems
\thanks{
\textbf{979-8-3503-4057-0/23/\$31.00 \copyright2023 IEEE.}
This work has been partly supported by the research project ``BotDozer,'' which is funded by EPSRC and EquipmentShare.
}
}

\author{\IEEEauthorblockN{Jie Lin}
\IEEEauthorblockA{\textit{School of EEECS} \\
\textit{Queen's University Belfast}\\
Belfast, Northern Ireland, UK \\
jlin20@qub.ac.uk}
\and
\IEEEauthorblockN{Ruairi Moran}
\IEEEauthorblockA{\textit{School of EEECS  \& i-AMS} \\
\textit{Queen's University Belfast}\\
Belfast, Northern Ireland, UK \\
rmoran05@qub.ac.uk}
\and
\IEEEauthorblockN{Pantelis Sopasakis}
\IEEEauthorblockA{\textit{School of EEECS \& i-AMS} \\
\textit{Queen's University Belfast}\\
Belfast, Northern Ireland, UK \\
p.sopasakis@qub.ac.uk}
}

\maketitle

\begin{abstract}
Nonlinear (deterministic and stochastic) optimal control problems are often solved on embedded devices using first-order numerical optimisation methods. The gradient computation accounts for a significant part of the computation cost per iteration; this is often performed with reverse-mode automatic differentiation and software libraries such as \casadi{} can be used to generate C code for this computation. In this paper, we propose a simple \textit{ad hoc} and highly parallelisable algorithm for the computation of the gradient of the total cost for deterministic and stochastic scenario-based optimal control problems. We also present \gradgen: an open-source Python package that generates Rust code for the gradient computation. The proposed method leads to a faster performance compared to \casadi{} and a significant reduction in generated code.
\end{abstract}


\section{Introduction}

\subsection{Background and motivation}
Gradient-based methods --- such as line search methods and various quasi-Newtonian methods \cite{nocedal_wright::2006} --- are popular for the solution of nonlinear optimal control problems (OCPs), especially on embedded devices since the involved steps are simple (typically no linear systems need to be solved), and the memory requirements are low. Then, the computation of the gradient of the cost function accounts for a major part of the computational cost at each iteration of the method. 

The gradient can be computed numerically, but that would introduce errors that could compromise the convergence of the numerical optimisation method. Gradients can also be computed using symbolic algebra methods, but symbolic differentiation often leads to overly lengthy expressions, which are inefficient for applications of numerical optimisation. 

Automatic differentiation (AD) is a class of methods for the algorithmic computation of the gradient \cite{margossian::2019}. There are two main AD approaches: forward-mode --- which is, in essence, equivalent to symbolic differentiation \cite{Laue2019OnTE} --- and reverse-mode AD. In forward-mode AD, the computational graph of a function $f:\R^n\to\R$ is traversed forwards, from the inputs $x_1, \ldots, x_n$ to the output $f(x)$, while the partial gradients $\partial f/\partial x_i$, $i=1,\ldots, n$, are computed. In reverse-mode AD, the computation that defines the function is traversed backwards, from $f(x)$ to $x$, while the chain rule is applied. Reverse-mode AD is more efficient compared to forward-mode \cite{giles::2008}.

Several open-source AD
implementations are available, such as JAX \cite{jax::2018} 
and within TensorFlow \cite{tensorflow::2015}.
\casadi{} \cite{andersson+::2019} is a Python package that can generate library-free C code that has been widely used in embedded applications such as 
autonomous ground vehicles \cite{agv::2018,Silano::2018}, 
aerial vehicles \cite{sina:2020,mav::2019}, and 
quadruped robots \cite{cheetah::2019}
to name a few. \casadi{} is also used to compute gradients in numerical optimisation software such as \texttt{OpEn} \cite{open::2020} and \texttt{acados} \cite{verschueren::2021}. However, the C code auto-generated by \casadi{} can often be tens to hundreds of thousands of lines of code long. 

The efficient computation of the cost gradient becomes even more important in problems of very large dimension such as stochastic scenario-based optimal control where the number of control actions increases exponentially with the prediction horizon \cite{sampathirao::2021}. The structure of scenario-based problems has been exploited to devise parallelisable numerical algorithms \cite{sampathirao::2021,sampathirao::2018} for convex formulations, but to the best of the authors' knowledge such developments have not been generalised to nonlinear formulations.

In this paper, we propose reverse-mode-type AD algorithms in the deterministic and stochastic cases, that allow fast and low-memory computation of the gradient of an OCP's cost function with respect to the control actions. We show that the complexity of the proposed algorithms scales gracefully with the problem size. In the stochastic case, the proposed algorithm can be massively parallelised. We also present an open-source code generator that generates fast memory-safe Rust code \cite{matsakis::2014} that is suitable for embedded applications.


\subsection{Notation}
Let $\N_{[k_1, k_2]}$ denote the integers in $[k_1, k_2]$.
We denote the transpose of a matrix $A$ by $A^{\tr}$.
Let $g: X \to Y$ and $f: Y \to Z$. Then we define the composition $f \circ g:X \to Z$ with $(f \circ g)(x) = f(g(x))$.
For $x\in\R^n$, note that $\|x\|$ denotes the Euclidean norm, that is, $\|x\| = \sqrt{x^{\tr} x}$.

\subsection{The Jacobian matrix and the chain rule}
Consider a function $f:\R^n\ni x \mapsto (f_1(x), \ldots, f_m(x))\in\R^m$.
The Jacobian matrix of a function $f$ is defined as the function $Jf:\R^n\to\R^{m\times n}$ with 
\begin{equation}
    Jf(x) = \begin{bmatrix}
    \frac{\partial f_1(x)}{\partial x_1} & \cdots &\frac{\partial f_1(x)}{\partial x_n}
    \\
    \vdots & \ddots & \vdots
    \\
    \frac{\partial f_m(x)}{\partial x_1} & \cdots &\frac{\partial f_m(x)}{\partial x_n}
    \end{bmatrix},
\end{equation}
provided that all partial derivatives exist. For functions $f:\R^n\to\R$ the gradient of $f$ is defined as the function $\nabla f:\R^n\to\R^n$ defined by $\nabla f(x) = Jf(x)^\tr$.

Suppose that $f:\R^p\to\R^m$ and $g:\R^n\to\R^p$ are differentiable. Then, the Jacobian matrix of the composition $f\circ g:\R^n\to\R^m$ is given by 
\begin{equation}
    J(f\circ g)(x) = Jf(g(x))Jg(x).
\end{equation}
This is the chain rule for Jacobian matrices. 
The notation $Jf(g(x))$ means the Jacobian of $f$ evaluated at $g(x)$.

Let $h_1: \R^n \ni x \mapsto s_1 \in \R^{n_1}$, and $h_2: \R^n \ni x \mapsto s_2 \in \R^{n_2}$. Define $g: \R^{n_1} \times \R^{n_2} \ni (s_1, s_2) \mapsto g(s_1, s_2) \in \R^m$. Then, $g$ is a function composed of two functions, and so we define $f(x) = g \circ (h_1(x), h_2(x))$. The Jacobian of $f$ at $x$ is
\begin{multline}    
    Jf(x) = J_{s_1} g(h_1(x), h_2(x)) Jh_1(x) 
    \\
    + J_{s_2} g(h_1(x), h_2(x)) Jh_2(x).
    \label{eq:chain_two_args}
\end{multline}

\section{Deterministic OCPs}\label{sec:deterministic}
Consider a --- generally nonlinear --- OCP of the general form 
\begin{subequations}
\begin{align}
    \mathbb{P}_N(x){}:{}
    \minimise_{
        \substack{
        u_0, \ldots, u_{N-1}\\
        x_0, \ldots, x_{N}
        }
    }& \sum_{t=0}^{N-1} \ell(x_t, u_t) + V_f(x_N),
    \\
    \subjto\ {}&x_{t+1} = f(x_t, u_t), t\in\N_{[0, N-1]},
    \label{eq:ocp-deterministic:dynamics}
    \\
    &x_0 = x,
\end{align}
\end{subequations}
where $x_t\in\R^{n_x}$ is the system state, 
$u_t\in\R^{n_u}$ is the input, 
and the system dynamics $f:\R^{n_x}\times\R^{n_u}\to\R^{n_x}$, 
the stage cost function $\ell:\R^{n_x}\times\R^{n_u}\to\R^{n_x}$, and the terminal cost function $V_f:\R^{n_x}\to\R^{n_x}$ are differentiable.

We can cast this problem as an unconstrained optimisation problem by eliminating the sequence of states. This is referred to as the single shooting formulation \cite{andersson+::2019}. We define $u = (u_0, u_1, \ldots, u_{N-1})\in\R^{Nn_u}$, $F_0(x, u) = x$, and recursively 
\begin{equation}
    F_{t+1}(x, u) = f(F_t(x, u), u_{t}).
\end{equation}
We see that $x_{t} = F_t(x, u)$ is the state at time $t$ starting from the initial state $x$ upon the action of the sequence $u$. This allows us to define the total cost function $V_N:\R^{n_x}\times\R^{Nn_u}\to\R$ as 
\begin{equation}
    V_N(x, u) {}={} \sum_{t=0}^{N-1} \ell(F_t(x, u), u_t) + V_f(F_N(x, u)),
\end{equation}
and $\mathbb{P}_N$ can be written as 
\begin{equation}
    \mathbb{P}_N(x)
{}:{}
    \minimise_{u\in\R^{Nn_u}}V_N(x, u).
\end{equation}
In embedded optimisation we typically aim at determining a stationary point of this problem, that is, a sequence of control actions $u^\star(x)$ such that $\nabla V_N(x, u^\star(x)) = 0$. To that end, any first-order optimisation method can be used; e.g., a line search gradient method, a quasi-Newtonian method \cite{nocedal_wright::2006}, or the bespoke PANOC method for OCPs \cite{panoc::2017}. In all these methods, the most computationally demanding part is typically the gradient of the total cost function.

A straightforward application of the chain rule for the computation of the gradient of $V_N$ is possible, but would necessitate the computation and storage of several matrices as well as several matrix-matrix/vector products. Instead, we will exploit the structure of the problem and will work in a fashion similar to the reverse-mode AD method \cite{giles::2008}. 

First, let us introduce some convenient notation. Given a vector $d\in\R^{n_x}$ we define $f^{x}_{t}(d) = J_x f(x_t, u_t)^\tr d$ and $f^{u}_{t}(d) = J_u f(x_t, u_t)^\tr d$, for $t\in\N_{[0, N-1]}$, where $x_t = F_t(x, u)$. We also define $\ell^{x}_{t} = \nabla_x \ell(x_t, u_t)$ and $\ell^{u}_{t} = \nabla_u \ell(x_t, u_t)$, for $t\in\N_{[0, N-1]}$, and $V_f' = \nabla V_f(x_N)$.

We start by determining the gradient of $V_N$ with respect to $u_{N-1}$, which, by \eqref{eq:chain_two_args}, is 
\begin{align}
    &\nabla_{u_{N-1}}V_N(x, u)
    \notag\\
    &\ {}={} \nabla_{u_{N-1}}[\ell(x_{N-1}, u_{N-1}) + V_f(f(x_{N-1}, u_{N-1}))]
    \notag\\
    &\ {}={} \ell^{u}_{N-1} + f^{u}_{N-1}(V_f').
\end{align}
Likewise, we can determine $\nabla_{u_{N-2}}V_N$. It is 
\begin{equation}
     \nabla_{u_{N-2}}V_N(x, u) 
     {}={}
     \ell^{u}_{N-2}
     {}+{}
     f^{u}_{N-2} (\ell^{x}_{N-1}
     {}+{}
      f^{x}_{N-1} (V_f')).
      \label{eq:diff_u_N_2}
\end{equation}
Note, however, that $V_f'$ has already been computed. Moving on to stage $N-3$, we have 
\begin{multline}
    \nabla_{u_{N-3}}V_N(x, u) 
    {}={}
    \\
    \ell^{u}_{N-3}
    {}+{}
    f^{u}_{N-3}(
    \ell^{x}_{N-2} 
    + 
    f^{x}_{N-2}(
    \ell^{x}_{N-1} 
    + 
    f^{x}_{N-1} (V_f'))).
\end{multline}
Note that the quantity $\ell^{x}_{N-1} + f^{x}_{N-1}(V_f')$ has already been computed in \eqref{eq:diff_u_N_2}. A pattern is emerging, which we state in the following proposition.
\begin{proposition}
Let $a_{N-1} = V_f'$. Then,
\begin{subequations}
    \begin{align}
        \nabla_{u_{t}} V_N(x, u) {}={}& \ell^u_t + f^{u}_t(a_t),
        \\
        a_{t-1} {}={}& \ell^x_t + f^x_t(a_t),
    \end{align}
\end{subequations}
for $t=N-1, N-2, \ldots, 0$.
\end{proposition}
This leads to Alg.~\ref{alg:gradient-deterministic} for the computation of the gradient of the total cost function, which is a special case of \cite[Alg. 1]{panoc::2017}.

\begin{algorithm}[htpb]
 \caption{Computation of the gradient of $V_N$ with respect to $u$ in the deterministic case}\label{alg:gradient-deterministic}
 \begin{algorithmic}[1]
 \renewcommand{\algorithmicrequire}{\textbf{Input:}}
 \renewcommand{\algorithmicensure}{\textbf{Output:}}
 \REQUIRE Initial state, $x$, and sequence of control actions, $u$
 \ENSURE  $\nabla V_N(x, u)$
  \FOR {$t = 0$ to $N-1$}
  \STATE $x_{t+1} \gets f(x_t, u_t)$
  \ENDFOR
  \STATE $a_{N-1} {}\gets{} V_f'$
  \FOR {$t = N-1$ to $1$}
    \STATE  $\nabla_{u_{t}} V_N(x, u) {}\gets{} \ell^u_t + f^{u}_t(a_t)$
    \STATE  $a_{t-1} {}\gets{} \ell^x_t + f^x_t(a_t)$
  \ENDFOR
  \STATE  $\nabla_{u_{0}} V_N(x, u) {}\gets{} \ell^u_0 + f^{u}_0(a_0)$
 \RETURN $\nabla V_N(x, u) = (\nabla_{u_{t}} V_N(x, u))_{t=0,\ldots, N-1}$
 \end{algorithmic}
 \end{algorithm}

 The variables $a_t$ in Alg.~\ref{alg:gradient-deterministic} correspond to the adjoint variables in reverse-mode AD. Note that Alg.~\ref{alg:gradient-deterministic} requires exactly $N$ invocations of $f$ and the Jacobians of $f$ and $\ell$ with respect to $x$ and $u$, that is, it scales linearly with $N$.

\section{Stochastic scenario-based OCPs}\label{sec:stochastic}
In this section, we turn our attention to stochastic scenario-based OCPs \cite{shapirolectures::2014}. We consider a dynamical system of the form 
\begin{equation}
    x_{t+1} = f(x_t, u_t, w_t),\label{eq:system_dynamics}
\end{equation}
for $t\in\N_{[0, N-1]}$, where $w_t$ is a finitely-supported random process acting on the system as a disturbance. At every time $t$, having measured $x_t$, we apply a control action $u_t = u_t(x_t)$. This entails a cost $\ell(x_t, u_t, w_t)$, which is a random variable. At the end of the horizon --- at $t=N$ --- a terminal cost, $V_f(x_N)$, is incurred.

The evolution of the system states and costs can be captured by the structure of a scenario tree as shown in \figurename~\ref{fig:scenario_tree}. The possible realisations of the state of the system correspond to the nodes of the scenario tree, while the edges correspond to the realisations of the disturbance, $w_t$.

\tikzset{main node/.style={circle,fill=blue!20,draw,minimum size=1cm,inner sep=0pt},
}

\begin{figure}[htpb!]
    \centering
    \begin{tikzpicture}[scale=0.55]
    \draw[thick, dashed, gray!60] (0,-7.3) -- (0,5);
    \draw[thick, dashed, gray!60] (5,-7.3) -- (5,5);
    \draw[thick, dashed, gray!60] (10,-7.3) -- (10,5);
    \node (node0) at (0,0)   [Nodes] {0};
    \node (node1) at (5,3)   [Nodes] {1};
    \node (node2) at (5,-3)  [Nodes] {2};
    \node (node3) at (10,4)  [Nodes] {3};
    \node (node4) at (10,2)  [Nodes] {4};
    \node (node5) at (10,0)  [Nodes] {5};
    \node (node6) at (10,-2) [Nodes] {6};
    \node (node7) at (10,-4) [Nodes] {7};
    \node (node8) at (10,-6) [Nodes] {8};
    \draw[arrow] (node0) to (node1);
    \draw[arrow] (node0) to (node2);
    \draw[arrow] (node1) to (node3);
    \draw[arrow] (node1) to (node4);
    \draw[arrow] (node2) to (node5);
    \draw[arrow] (node2) to (node6);
    \draw[arrow] (node2) to (node7);
    \draw[arrow] (node2) to (node8);
    \draw[axis] (-1,-7.3) -- (11,-7.3);
    \node (t0) at (0,-7.8) {\color{gray}$t=0$};
    \node (t1) at (5,-7.8) {\color{gray}$t=1$};
    \node (t2) at (10,-7.8) {\color{gray}$t=2$};
    \node [left=0 of node0] (x^0) {$x^0$};
    \node [left=0 of node1] (x^1) {$x^1$};
    \node [left=0 of node2] (x^2) {$x^2$};
    \node [right=1cm of node0] (u^0) {$u^0$};
    \draw[arrow] (node0) to (u^0);
    \node [right=1cm of node1] (u^1) {$u^1$};
    \draw[arrow] (node1) to (u^1);
    \node [right=1cm of node2] (u^2) {$u^2$};
    \draw[arrow] (node2) to (u^2);
    \node [yshift=2.5cm,right=-0.3cm of node0] (w_0) {$w_{0}$};
    \node [yshift=2.5cm,right=2.4cm of node0] (w_1) {$w_{1}$};
    \node [rotate=30]  (ellx^0u^0w_0) at (2.7,2.1)  {\color{ForestGreen}$\ell^1$};
    \node [rotate=-30] (ellx^0u^0w_0) at (2.7,-2.1) {\color{ForestGreen}$\ell^2$};
    \node [rotate=11]  (ellx^1u^1w_1) at (7.8,4.1)  {\small\color{ForestGreen}$\ell^3$};
    \node [rotate=-11] (ellx^1u^1w_1) at (7.8,2)  {\small\color{ForestGreen}$\ell^4$};
    \node [rotate=30]  (ellx^2u^2w_1) at (7.8,-0.8) {\small\color{ForestGreen}$\ell^5$};
    \node [rotate=12]  (ellx^2u^2w_1) at (7.9,-2.0) {\small\color{ForestGreen}$\ell^6$};
    \node [rotate=-12] (ellx^2u^2w_1) at (7.9,-4.0) {\small\color{ForestGreen}$\ell^7$};
    \node [rotate=-30] (ellx^2u^2w_1) at (7.8,-5.2) {\small\color{ForestGreen}$\ell^8$};
    \node [,right=0 of node3] (ellx^3) {\color{ForestGreen}$V_f^3$};
    \node [right=0 of node4] (ellx^4) {\color{ForestGreen}$V_f^4$};
    \node [right=0 of node5] (ellx^5) {\color{ForestGreen}$V_f^5$};
    \node [right=0 of node6] (ellx^6) {\color{ForestGreen}$V_f^6$};
    \node [right=0 of node7] (ellx^7) {\color{ForestGreen}$V_f^7$};
    \node [right=0 of node8] (ellx^8) {\color{ForestGreen}$V_f^8$};
    \node [yshift=-0.5cm, left=-0.3 of node0] (pi^0) {$\pi^0$};
    \node [yshift=-0.5cm, left=-0.3 of node1] (pi^1) {$\pi^1$};
    \node [yshift=-0.5cm, left=-0.3 of node2] (pi^2) {$\pi^2$};
    \node [yshift=-0.5cm, left=-0.3 of node3] (pi^3) {$\pi^3$};
    \node [yshift=-0.5cm, left=-0.3 of node4] (pi^4) {$\pi^4$};
    \node [yshift=-0.5cm, left=-0.3 of node5] (pi^5) {$\pi^5$};
    \node [yshift=-0.5cm, left=-0.3 of node6] (pi^6) {$\pi^6$};
    \node [yshift=-0.5cm, left=-0.3 of node7] (pi^7) {$\pi^7$};
    \node [yshift=-0.5cm, left=-0.3 of node8] (pi^8) {$\pi^8$};
\end{tikzpicture}
    \caption{Scenario tree with $N=2$. Regarding the system dynamics, note that $x^1=f(x^0, u^0, w^1)$ and $x^2 = f(x^0, u^0, w^2)$. Likewise, for the stage costs note that $\ell^2 = \ell(x^0, u^0, w^2)$. For the terminal costs, note that $V_{f}^{3} = V_{f}(x^{3})$.}
    \label{fig:scenario_tree}
\end{figure}
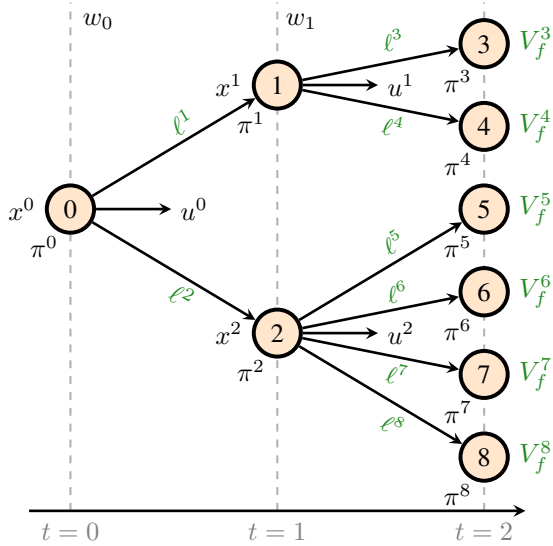

The scenarios are enumerated with an index $i$, where $i=0$ corresponds to the unique node at stage $t=0$ which is referred to as the \textit{root} node. 
The set of nodes at stage $t\in\N_{[0, N]}$ are denoted by $\nodes(t)$, where $\nodes(0)=\{0\}$. 
The nodes at stage $t=N$ are referred to as the \textit{leaf} nodes of the tree, and the nodes that are not at this stage are referred to as the \textit{nonleaf} nodes.
We also denote $\nodes(t_1, t_2) = \bigcup_{\tau=t_1}^{t_2}\nodes(\tau)$.
Each node $i\in\nodes(1, N)$ is connected to a single node at stage $t-1$ known as the ancestor of $i$, denoted by $\anc(i)$. Conversely, every node $i\in\nodes(0, N-1)$ at a stage $t$ leads up to a set of nodes $\child(i) \subseteq \nodes(t+1)$. 
The probability of occurrence of node $i$ is denoted by $\pi^i > 0$. Lastly, every node of the tree is associated with a state $x^i$, every nonleaf node with an input $u^i$, and every edge --- so, all but the root node --- with a disturbance $w^i$. Note that $w^i$ is the event that happens at node $\anc(i)$ that causes node $i$ to happen.

The dynamics \eqref{eq:system_dynamics} on the scenario tree can be written
\begin{equation}
    x^{i_+} = f(x^i, u^i, w^{i_+}),
\end{equation}
for $i\in\nodes(0, N-1)$ and $i_+ \in \child(i)$. Equivalently,
\begin{equation}
    x^{i} = f(x^{\anc(i)}, u^{\anc(i)}, w^i),
\end{equation}
for $i\in\nodes(1, N)$.

We can now formulate the stochastic scenario-based OCP 
\begin{subequations}
    \begin{align}
        \minimise_{
        \substack{
        (u^i)_{i\in\nodes(0, N-1)},\\
        (x^i)_{i\in\nodes(0, N)}
        }
    }&\ \sum_{i\in\nodes(1,N-1)} \pi^i \ell(x^{\anc(i)}, u^{\anc(i)}, w^i) 
    \notag\\
    &\quad+ \sum_{j\in\nodes(N)}\pi^j V_f(x^j),
    \\
    \subjto&\  x^{i} = f(x^{\anc(i)}, u^{\anc(i)}, w^i),
    \notag\\
    &\qquad i\in\nodes(1, N),
    \\
    &x^0 = x.
    \end{align}
\end{subequations}
Hereafter we assume that $f$ and $\ell$ are differentiable in $(x, u)$ and $V_f$ is differentiable.

Let $u = (u^i)_{i\in\nodes(0, N-1)}$. We define $V_N(x, u)$ to be the total cost function in the above problem. We will determine the gradient of $V_N$ with respect to all $u^i$. At this point, let us introduce the notation $\ell_u^{i} = \nabla_u \ell(x^{\anc(i)}, u^{\anc(i)}, w^i)$ and $\ell_x^{i} = \nabla_x \ell(x^{\anc(i)}, u^{\anc(i)}, w^i)$, for $i\in\nodes(1, N)$.
Moreover, for $d\in\R^{n_x}$ we define $f_x^i(d) = J_x f(x^{\anc(i)}, u^{\anc(i)}, w^i)^\tr d$ and similarly, $f_u^i(d)$, for $i\in\nodes(1, N)$. Lastly, we define  $V_f^{\prime i} = \nabla V_f(x^i)$, for $i\in\nodes(N)$.

Firstly, for all $k_0\in\nodes(N)$ we compute $a^{k_0} = V_{f}^{\prime k_0}$. Then for $k_1\in\nodes(N-1)$, we compute
\begin{subequations}\label{eq:stochastic:n-1}
\begin{equation}\label{eq:nabla-u-k1}
    \nabla_{u^{k_1}}V_N
    {}={} 
    \sum_{k_0\in\child(k_1)} \pi^{k_0} (\ell_{u}^{k_0} + f_{u}^{k_0} (
    \underbracket[0.5pt]{
    V_{f}^{\prime k_0}}_{a^{k_0}})),
\end{equation}
and
\begin{equation}
    a^{k_1} {}={} 
    \sum_{k_0\in\child(k_1)}
    \pi^{k_0}(\ell_x^{k_0} 
    + f_x^{k_0}(a^{k_0})).
\end{equation}
\end{subequations}
Note that $a^{k_1}$ is similar to the gradient in \eqref{eq:nabla-u-k1}, except with respect to $x$ instead of $u$.
Next, let us take $k_2\in\nodes(N-2)$. We compute
\begin{align}\label{eq:nabla-u-k2}
    \nabla_{u^{k_2}}V_N
    {}={}&
    \sum_{k_1\in\child(k_2)}\hspace{-.5em}
    \pi^{k_1}\ell_u^{k_1}
    \notag
    \\
    & \quad {}+{}
    \sum_{\substack{k_1\in\child(k_2)\\k_0\in\child(k_1)}}
    \hspace{-.5em}
    \pi^{k_0}(f_u^{k_1}(\ell_x^{k_0}) 
    + f_u^{k_1}(f_x^{k_0}(V_f^{\prime k_0})))
    \notag
    \\
    {}={}&
    \sum_{k_1\in\child(k_2)}\hspace{-.5em}
    \pi^{k_1}\ell_u^{k_1}
    \notag
    \\
    & \quad {}+{}
    \sum_{k_1\in\child(k_2)}\hspace{-.5em}
    f_u^{k_1}
    \Big(\hspace{-1.1em}
    \underbracket[0.5pt]{
    \sum_{k_0\in\child(k_1)}\hspace{-.5em}
    \pi^{k_0}(\ell_x^{k_0} 
    + f_x^{k_0}(V_f^{\prime k_0}))}_{a^{k_1}}
    \hspace{-.3em}\Big)
    \notag
    \\
    {}={}&
    \sum_{k_1\in\child(k_2)}\hspace{-.5em}
    \pi^{k_1}\ell_u^{k_1}
    {}+{}
    f_u^{k_1} (a^{k_1}),
\end{align}
and
\begin{align}
    a^{k_2} {}={} \sum_{k_1\in\child(k_2)}
    \pi^{k_1}\ell_x^{k_1} 
    + f_x^{k_1}(a^{k_1}).
\end{align}
Again, note that $a^{k_2}$ is similar to the gradient in \eqref{eq:nabla-u-k2}, except with respect to $x$ instead of $u$. For $\nabla_{u^i}V_N$ with $i\in\nodes(N-1)$, we use \eqref{eq:stochastic:n-1}.
Then, moving backwards, we can show recursively that for $i\in\nodes(0, t)$ with $t\in\N_{[0, N-2]}$
\begin{subequations}
\begin{equation}
    \nabla_{u^{i}}V_N 
    {}={} 
    \sum_{i_+\in\child(i)}
    \pi^{i_+}\ell_u^{i_+} 
    + f_u^{i_+}(a^{i_+}),
\end{equation}
where the adjoint vectors $a^{i}$, for $i\in\nodes(1, N-2)$ are given by
\begin{equation}
    a^{i} 
    {}={} 
    \sum_{i_+\in\child(i)}
    \pi^{i_+}\ell_x^{i_+} 
    + f_x^{i_+}(a^{i_+}).
\end{equation}
\end{subequations}

Many of the involved computations can be parallelised. This is evident in Alg.~\ref{alg:gradient-stochastic} where the parallelisation is represented using the fork-join model \cite[Chapter 2]{pacheco:2011}. In fact, the computation of the gradient of the total cost lends itself to a lockstep implementation on a GPU or FPGA.

\begin{algorithm}[htpb!]
 \caption{Parallelisable computation of the gradient of $V_N$ with respect to $u$ in the stochastic case}\label{alg:gradient-stochastic}
 \begin{algorithmic}[1]
 \renewcommand{\algorithmicrequire}{\textbf{Input:}}
 \renewcommand{\algorithmicensure}{\textbf{Output:}}
 \REQUIRE Initial state, $x$, and sequence of control actions, $u$
 \ENSURE  $\nabla V_N(x, u)$
  \FOR {$i\in\nodes(1, N)$}
  \STATE $x^i \gets f(x^{\anc(i)}, u^{\anc(i)}, w^i)$
  \ENDFOR
  \FOR {$i\in\nodes(N)$ (\textit{in parallel})}
    \STATE \textbf{fork} $a^{i} \gets V_f^{\prime i}$
  \ENDFOR
  \STATE \textbf{join}
  \FOR {$i\in\nodes(N-1)$ (\textit{in parallel})}
    \STATE \textbf{fork:} 
    \STATE \quad $\nabla_{u^i}V_N {}\gets{}
        \sum_{i_+\in\child(i)}\pi^{i_+}(\ell_u^{i_+} 
        + f_u^{i_+}a^{i_+})$
        \STATE \quad $a^{i} {}\gets{} \sum_{i_+\in\child(i)}
            \pi^{i_+}(\ell_x^{i_+} + f_x^{i_+}a^{i_+})$   
  \ENDFOR
  \STATE \textbf{join}
  \FOR {$t = N-2$ \TO $1$}
    \FOR {$i \in \nodes(t)$ (\textit{in parallel})}
        \STATE \textbf{fork:} 
        \STATE \quad $\nabla_{u^i}V_N {}\gets{}
            \sum_{i_+\in\child(i)}\pi^{i_+}\ell_u^{i_+} 
            + f_u^{i_+}a^{i_+}$
        \STATE \quad $a^{i} {}\gets{} \sum_{i_+\in\child(i)}
            \pi^{i_+}\ell_x^{i_+} + f_x^{i_+}a^{i_+}$
    \ENDFOR
    \STATE \textbf{join}
  \ENDFOR
  \STATE  $\nabla_{u^{0}}V_N {}\gets{}
        \sum_{i_+\in\child(0)}
        \pi^{i_+}\ell_u^{i_+}{}+{}f_u^{i_+}a^{i_+}$
 \RETURN $\nabla V_N(x, u) = (\nabla_{u^i} V_N(x, u))_{i\in\nodes(0,N-1)}$
 \end{algorithmic}
 \end{algorithm}

 Note that the complexity of Alg.~\ref{alg:gradient-stochastic} is linear with respect to the total number of nodes of the tree.

\section{gradgen: Code generation for embedded computations}\label{sec:gradgen}

The open-source Python package \gradgen{}\footnote{The source code is available at \url{https://github.com/QUB-ASL/gradgen}.} is distributed under the MIT License and can be easily installed with \texttt{pip}.

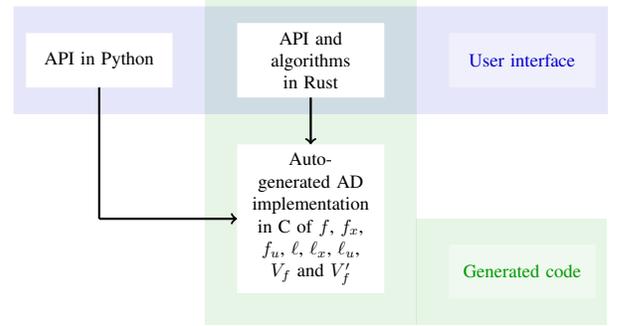
\begin{figure}[htbp!]
    \centering
    \resizebox{0.45\textwidth}{!}{
        \begin{tikzpicture}
    \fill [fill=mygreen,very thick,opacity=0.1] (-2,-5) rectangle (2,1.2);
    \fill [fill=mygreen,very thick,opacity=0.1] (2,-5) rectangle (5.6,-3);

    \draw [draw=myblue,fill=myblue,very thick,opacity=0.1] (-5.6,-1.0) rectangle (5.6,1.0);
    
    \node[block] (python) at (-4, 0) {API in Python};
    
    \node[block,very thick] (gradgen) at (0, 0) {API and algorithms in Rust};
    
    \node[block,very thick] (casadi) at (0, -3) {Auto-generated AD implementation in C of $f$, $f_{x}$, $f_{u}$, $\ell$, $\ell_{x}$, $\ell_{u}$, $V_{f}$ and $V_{f}'$};
    
    \node[block,text=myblue, opacity=0.4, text opacity=1] (ui) at (4, 0) {User interface};
    
    \node[block,text=mygreen, opacity=0.4, text opacity=1] (gen) at (4, -4) {Generated code};
    
    \node[] (dummyarrow) at (-4, -3) {};
    
    \draw[-, very thick] (python.south) -- (dummyarrow.center);
    \draw[->, very thick] (dummyarrow.center) -- (casadi.west);
    \draw[->, very thick] (gradgen) -- (casadi);
\end{tikzpicture}
    }
    \caption{Multi-layered structure of the generated code and Python API.}
    \label{fig:code_structure}
\end{figure}

The user can then import \gradgen{} to their Python project, where they specify the system dynamics and cost function in symbolic form. The Jacobians $\ell^x, \ell^u, f^x, f^u,$ and $V_f'$ are determined using \casadi. Then using Jinja2 templates, \texttt{gradgen} generates Rust code that implements Alg.~\ref{alg:gradient-deterministic} or \ref{alg:gradient-stochastic}.
The modular structure of the generated code is presented in Fig.~\ref{fig:code_structure}.
The top layer is the generated \gradgen{} API, comprising Rust implementations of Alg.~\ref{alg:gradient-deterministic} and \ref{alg:gradient-stochastic}. 



\section{Results}\label{sec:simulations}
Here we demonstrate how \gradgen{} performs on two systems: an inverted pendulum (two states, one input) and a ball-and-beam (four states, one input), in both deterministic and stochastic cases. All simulations are run on an Apple Macbook Pro with a \unit[2]{GHz} quad-core Intel core i5 processor and \unit[16]{GB} of RAM.
Below, the system state is denoted $x$ and the control action $u$. For both examples, the stage cost  is $\ell(x,u) = \|x\|^2 + \|u\|^2$, and the terminal cost  is $V_f(x) = 10\|x\|^2$.

\subsection{System \Romannum{1}: inverted pendulum}
\begin{figure}[htpb!]
    \centering
    \input{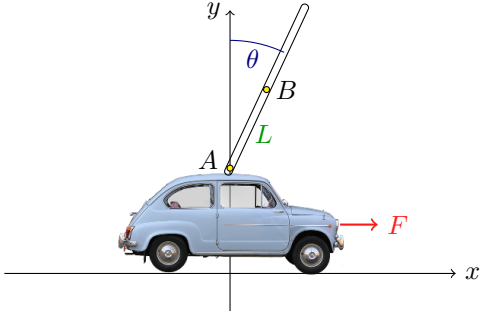}
    \caption{System \Romannum{1}: inverted pendulum ($n_x=2, n_u=1$). The mid-point of the rod is denoted $B$. The force applied to the car is denoted by $F$. The deviation angle, $\theta$, is the angle the rod makes with the $y$ axis.}
    \label{fig:Inverted pendulum}
\end{figure}

Consider the inverted pendulum system in \figurename~\ref{fig:Inverted pendulum}.
The car can move along the $x$ axis, while the rod is free to rotate about the fulcrum, $A$, which is fixed to the car. The deviation angle, $\theta$, is the angle the rod makes with the $y$ axis.
The car moves under the horizontal force, $F$, in an attempt to keep the rod upright, that is, $\theta = 0$.

The car has a mass $M$.
The rod has a mass $m$ and half-length $L$ and the rod is isotropic.
The friction of rotation of the rod is negligible. 

The dynamics of the inverted pendulum is \cite{sopasakis_book::2023}
\begin{subequations}
\begin{align}
  \dot{\theta} 
  {}={}& 
  \omega,
  \\
  \dot{\omega} 
  {}={}& 
  -3 \frac{
  \tfrac{1}{2} m L \omega^2 \sin(2 \theta)
  + F \cos\theta
  - M_{\rm tot} g \sin\theta
  }
  {
  (4 M_{\rm tot} - 3 m \cos^2 \theta)L
  }.\label{eq:Inverted pendulum}
\end{align}
\end{subequations}
where $M_{\rm tot} = M + m$.

We discretise this system using the Euler method with sampling time $T_s$ to obtain a discrete-time system in the form of \eqref{eq:ocp-deterministic:dynamics}. We use the values $m = \unit[1]{kg}$, $M = \unit[3]{kg}$, $g = \unitfrac[9.81]{m}{s^2}$, and $T_s = \unit[0.01]{s}$.


\subsection{System \Romannum{2}: ball-and-beam}
\begin{figure}[htpb!]
    \centering
    \includegraphics[width=0.75\linewidth]{\pathtofigs/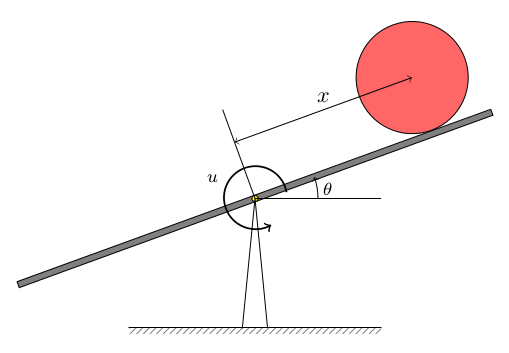}
    \caption{System \Romannum{2}: ball and beam ($n_x=4, n_u=1$). The position of the ball is denoted by $x$, the angle of the beam is $\theta$. The manipulated input is the torque $u$.}
    \label{fig:ball_beam_system}
\end{figure}

 Consider a ball of mass $m$ is placed on a beam which is poised on a fulcrum at its middle as shown in \figurename~\ref{fig:ball_beam_system}. 
We can control the system by applying a torque $u$ with respect to the fulcrum point.
The moment of inertia of the beam is denoted by $I$. The displacement $x$ of the ball from the 
midpoint can be measured with an optical sensor. 
The dynamical system is described by the following nonlinear differential equations
\begin{subequations}
\begin{align}
 \tfrac{7}{5}\ddot{x} + g\sin\theta - x\dot{\theta}^2 {}={}& 0,
 \\
 (mx^2 + I)\ddot{\theta} + 2mx\dot{x}\dot{\theta} + mgx\cos \theta {}={}& u.
\end{align}
\end{subequations}
This is a continuous-time system with state vector $(x, \dot{x}, \theta, \dot{\theta})$ and input $u$, which can be discretised with the Euler method and sampling period $T_s = \unit[0.01]{s}$.

\subsection{Stochastic systems}\label{sec:stochastic}
\begin{figure}[htpb!] 
    \centering
    \begin{tikzpicture}[scale=0.7]
    [->,>=stealth',shorten >=2pt, line width=0.7pt]
    \node [circle , draw] (zero) at (0, 0) {$\unit[10]{ms}$};
    \node [circle, draw] (one) at (2, 2) {$\unit[20]{ms}$};
    \node [circle, draw] (two) at (4, -1.8) {$\unit[100]{ms}$};

    \path[->] (zero) edge [bend left] node [above,yshift=2mm] {3\%} (one); 
    \path[->] (one) edge [bend left] node [below,xshift=3mm] {85\%} (zero);
    \draw[->] (zero) to[out=-10,in=170] node[above,xshift=4mm,yshift=-2mm] {2\%} (two);

    \path[->] (two) edge [bend left] node [below] {90\%} (zero); 
    \path[->] (one) edge [bend left] node [right] {5\%} (two);
    \draw[->] (two) to[out=120,in=280] node[right] {5\%} (one);

    \path (zero) edge [loop left] node {95\%} (zero);
    \path (one) edge [loop right] node {10\%} (one);
    \path (two) edge [loop right] node {5\%} (two);
\end{tikzpicture} 
    \caption{Markov chain of random sampling time. The sampling time $T_s =\unit[10]{ms}$ corresponds to normal operation, whereas $\unit[20]{ms}$ and $\unit[100]{ms}$ correspond to random time delays.} \label{fig:markov_chain}
\end{figure}
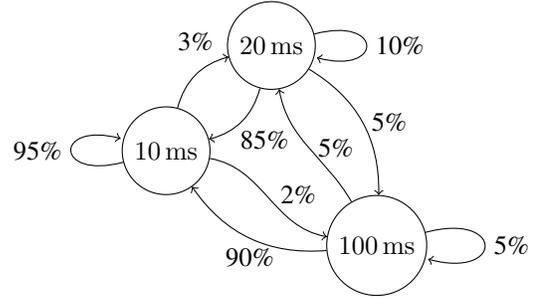

 Here we consider stochastic variants of the above systems by considering the case of random communication delays that render the sampling time, $T_s$, a random variable \cite{patrinos::2011}. We assume that $T_s$ can take the values $\unit[10]{ms}$, $\unit[20]{ms}$, and  $\unit[100]{ms}$, which are governed by the Markov chain shown in \figurename~\ref{fig:markov_chain} with uniform initial distribution.

\subsection{Benchmark results}
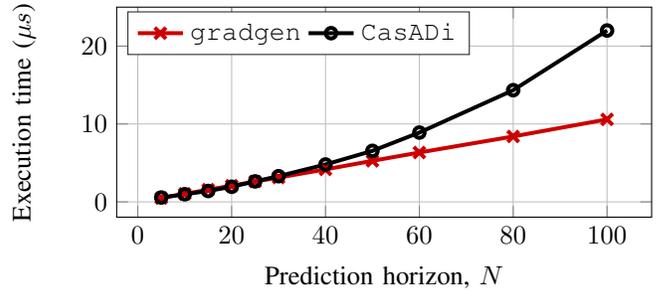
\begin{figure}[htpb!]
    \centering
    \begin{tikzpicture}%
\begin{axis}[%
    width=2.8in,
    height=1.1in,
    scale only axis,
    xminorticks=true,
    xlabel={Prediction horizon, $N$},
    ymax=25,
    yminorticks=true,
    ylabel={Execution time ($\mu{}s$)},
    xmajorgrids,
    xminorgrids,
    ymajorgrids,
    yminorgrids,
    legend columns=2, 
    legend style={
        at={(0.02,0.98)}, 
        anchor=north west, 
        legend cell align=left, 
        align=left, 
        draw=white!15!black}
]

\addplot [
    color=red!80!black, 
    line width=1.5pt, 
    mark=x,
    mark size=3pt, 
    mark options={solid, red!80!black}
    ]
  table[row sep=crcr]{%
    5.0000    0.49\\
   10.0000    1.03\\
   15.0000    1.57\\
   20.0000    2.08\\
   25.0000    2.62\\
   30.0000    3.12\\
   40 4.17\\
   50 5.27\\
   60 6.34\\ 
   80 8.40\\
   100 10.59\\
};
\addlegendentry{\gradgen}

\addplot [
    color=black, 
    line width=1.5pt, 
    mark=o,
    mark size=2pt, 
    mark options={solid, black}
    ]
  table[row sep=crcr]{%
    5	0.56\\ 
    10	0.98\\
    15	1.4\\ 
    20	1.94\\
    25	2.63\\
    30	3.30\\
    40 4.80\\
    50 6.55\\
    60 8.9\\
    80 14.34\\
    100 22\\
};
\addlegendentry{\casadi}

\end{axis}%

\end{tikzpicture}%
    \caption{Inverted pendulum example: average computation time for the gradient and comparison with \casadi{} (based on $10^6$ runs).}
    \label{fig:gradgen-vs-casadi-invpend}
\end{figure}

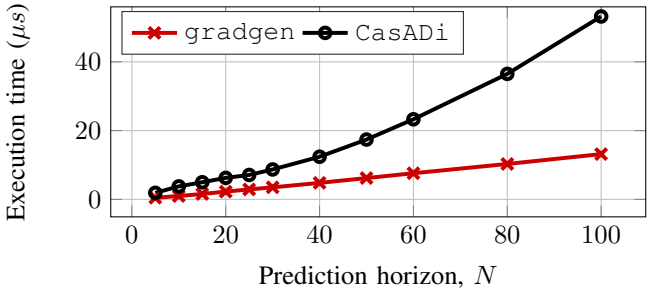
\begin{figure}[htpb!]
    \centering
    \begin{tikzpicture}%
\begin{axis}[%
    width=2.8in,
    height=1.1in,
    scale only axis,
    xminorticks=true,
    xlabel={Prediction horizon, $N$},
    ymax=56,
    yminorticks=true,
    ylabel={Execution time ($\mu{}s$)},
    xmajorgrids,
    xminorgrids,
    ymajorgrids,
    yminorgrids,
    legend columns=2, 
    legend style={
        at={(0.02,0.98)}, 
        anchor=north west, 
        legend cell align=left, 
        align=left, 
        draw=white!15!black}
]

\addplot [
    color=red!80!black, 
    line width=1.5pt, 
    mark=x,
    mark size=3pt, 
    mark options={solid, red!80!black}
    ]
  table[row sep=crcr]{%
    5.0000    0.47\\
   10.0000    0.98\\
   15.0000    1.58\\
   20.0000    2.25\\
   25.0000    2.9\\
   30.0000    3.5\\
   40.0000    4.8\\
   50.0000    6.2\\
   60   7.6\\
   80   10.3\\
   100  13.17\\
};
\addlegendentry{\gradgen}

\addplot [
    color=black, 
    line width=1.5pt, 
    mark=o,
    mark size=2pt, 
    mark options={solid, black}
    ]
  table[row sep=crcr]{%
    5.0000    1.9\\
   10.0000    3.78\\
   15.0000    4.98\\
   20.0000    6.25\\
   25.0000    7.1\\
   30.0000    8.7\\
   40.0000    12.4\\
   50.0000    17.4\\
   60 23.3\\
   80 36.5\\
   100 53.2\\
};
\addlegendentry{\casadi}

\end{axis}%

\end{tikzpicture}%
    \caption{Ball-and-beam example: average computation time for the gradient and comparison with \casadi{} (based on $10^6$ runs).}
    \label{fig:gradgen-vs-casadi-bnb}
\end{figure}

In Figures \ref{fig:gradgen-vs-casadi-invpend} and \ref{fig:gradgen-vs-casadi-bnb} we present the average execution time for the computation of the gradient using \gradgen's implementation of Alg.~\ref{alg:gradient-deterministic} and \casadi{} (in C) for the inverted pendulum and the ball-and-beam systems. In both figures we can observe that the execution time of \gradgen{} increases linearly with the prediction horizon. Indeed, the complexity of Alg.~\ref{alg:gradient-deterministic} is $\mathcal{O}(N)$. Instead, the execution time of \casadi{} seems to be growing superlinearly with the prediction horizon.

It should be noted that \casadi{} generates large C code files. For instance, for the ball-and-beam system with $N=60$, the \casadi{} auto-generated C file counts $46\unit{k}$ single lines of code (SLOCs). At $N=100$, it generates $125\unit{k}$ SLOCs and the compilation time exceeds $20\unit{mins}$. On the other hand, the total SLOCs of \gradgen{} (including header files) is $1188$, which is independent of the prediction horizon. For the inverted pendulum model, the auto-generated SLOCs by \casadi{} at $N=50$ and $N=100$ are $12.5\unit{k}$ and $44.9\unit{k}$, respectively.





\begin{figure}[htpb!]
    \centering
    \begin{tikzpicture}

\begin{axis}[%
    width=2.8in,
    height=1.0in,
    at={(1.011in,0.642in)},
    scale only axis,
    xmode=log,
    xmin=10,
    xmax=1250,
    xminorticks=true,
    xlabel={Number of nodes},
    ymode=log,
    ymin=1,
    ymax=400,
    yminorticks=true,
    ylabel={Execution time ($\mu{}s$)},
    axis background/.style={fill=white},
    xmajorgrids,
    xminorgrids,
    ymajorgrids,
    yminorgrids,
    legend style={
        at={(0.98,0.02)}, 
        anchor=south east, 
        legend cell align=left, 
        align=left, 
        draw=white!15!black}
]
\addplot [color=black, line width=1.5pt, mark=x,
mark size=3pt, mark options={solid, black}]
  table[row sep=crcr]{%
13  1.9\\
40	5.7\\
121 18\\
364 57\\
1093 179\\
};
\addlegendentry{Ball-and-beam}

\addplot [color=black, dashed, line width=1.5pt, mark=o, mark options={solid, black}]
  table[row sep=crcr]{%
13  1.6\\
40  5\\
121 16\\    
364 46\\
1093    147\\
};
\addlegendentry{Inv. pendulum}

\end{axis}

\end{tikzpicture}%
    \caption{Execution time for the computation of the gradient of the stochastic optimal control problem against the number of nodes.}
    \label{fig:stochastic-gradient-time}
\end{figure}
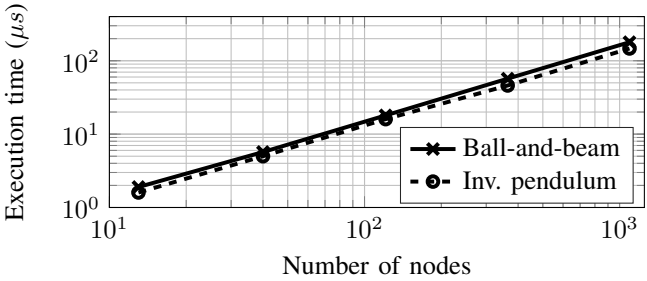

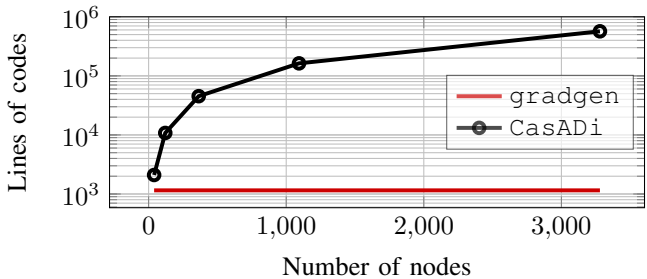
\begin{figure}[htpb!]
    \centering
    \begin{tikzpicture}%
\begin{semilogyaxis}[%
    width=2.8in,
    height=1.0in,
    scale only axis,
    xminorticks=true,
    xlabel={Number of nodes},
    ymax=1e6,
    yminorticks=true,
    ylabel={Lines of codes},
    xmajorgrids,
    xminorgrids,
    ymajorgrids,
    yminorgrids,
    legend columns=1, 
    legend style={ 
        anchor=east,
        at={(0.98,0.5)},
        legend cell align=left, 
        align=left, 
        draw=white!15!black,
        opacity=0.7,
        text opacity = 1}
]

\addplot [
    color=red!80!black, 
    line width=1.5pt, 
    ]
  table[row sep=crcr]{%
    40    1152\\
   3280    1152\\
};
\addlegendentry{\gradgen}

\addplot [
    color=black, 
    line width=1.5pt, 
    mark=o,
    mark size=2pt, 
    mark options={solid, black}
    ]
  table[row sep=crcr]{%
    40 2098\\
    121 10804\\
    364 45294\\
    1093 162910\\
    3280 567922\\
};
\addlegendentry{\casadi}

\end{semilogyaxis}%
\end{tikzpicture}%
    \caption{Ball-and-beam stochastic system: lines of code for the computation of the gradient of the stochastic optimal control problem against the number of nodes and comparison with \casadi.}
    \label{fig:SLOCs-vs-problem-size}
\end{figure}

In the stochastic case (Sec. \ref{sec:stochastic}), we compute the gradient of the total cost function, $\nabla V_N$, for $N=3,\ldots, 7$, using \gradgen's implementation of Alg.~\ref{alg:gradient-stochastic}. The computation times are shown in \figurename~\ref{fig:stochastic-gradient-time} where we see a linear increase with the number of nodes. \figurename~\ref{fig:SLOCs-vs-problem-size} shows the number of SLOCs of \casadi's auto-generated C file.

\vspace{-1em}

\vspace{-1em}

\section{Conclusions and future work}
In this paper we derived an algorithm for computing the gradient of a cost function with respect to a sequence of control actions in the deterministic case. We also derived a parallelisable algorithm for computing the gradient in the stochastic case. We simulated two systems that demonstrate the implementations of the deterministic and stochastic algorithms in the new open-source package \gradgen{}, which is available at \url{https://github.com/QUB-ASL/gradgen}. The deterministic algorithm is benchmarked against \casadi{}. The simulation results indicate that \gradgen{} can outperform \casadi{}, especially at large prediction horizons, \gradgen{} generates much smaller files than \casadi{}, and that the run time of \gradgen{} scales well with the problem size.
Future work will focus on implementing the stochastic algorithm in parallel on a GPU.

\bibliographystyle{ieeetr} 
\bibliography{bibdb}

\end{document}